\documentclass{amsart}
\usepackage[english]{babel}
\usepackage{graphicx}
\usepackage{amsmath}
\usepackage{amsthm}
\usepackage{amssymb}
\newtheorem{Main}{Theorem}[section]
\newtheorem{Lemma 1}[Main]{Lemma}
\newtheorem{Lemma 2}[Main]{Lemma}
\newtheorem{Lemma 3}[Main]{Lemma}
\newtheorem{Lemma 4}[Main]{Lemma}
\newtheorem{Lemma 5}[Main]{Lemma}
\begin{document}
\title{On $rth$ coefficient of divisors of $x^n-1$}
\address{School of Mathematics, Tata Institute of Fundamental Research Mumbai, India 400005}
\author{Sai Teja Somu}
\date{October 27, 2017}
\maketitle
\begin{abstract}
	Let $r,n$ be two natural numbers and let $H(r,n)$ denote the maximal absolute value of $r$th coefficient of divisors of $x^n-1$. In this paper, we show that $\sum_{n\leq x}H(r,n)$ is asymptotically equal to $c(r)x(\log x)^{2^r-1}$ for some constant $c(r)>0$. Furthermore, we give an explicit expression of $c(r)$ in terms of $r$. 
\end{abstract}
\section{Introduction}
Cyclotomic polynomials are irreducible divisors of $x^n-1$. The factorization of $x^n-1$ is as follows,
\[ x^n-1=\prod_{d|n}\Phi_d(x).\] Applying Mobius inversion formula we get $\Phi_n(x)=\prod_{d|n}(x^d-1)^{\mu(\frac{n}{d})}.$
Coefficients of $n$th cyclotomic polynomial has been the subject of study in \cite{Bachman},\cite{Bateman} and \cite{Bloom}.

In \cite{Pomerance}, the study of coefficients of divisors of $x^n-1$ has started. Recently, in \cite{Somu1} and \cite{Somu2} we proved that for every finite sequence of integers $n_1,\cdots,n_r$ there exists a natural number $n$ and a divisor $d(x)$ of $x^n-1$ of the form $d(x)=1+n_1x+\cdots +n_rx^r+\sum_{j=r+1}^{deg(d)}c_jx^j$.

Let $(f)_r$ denote $rth$ coefficient of the $f(x)$ for any formal power series $f(x)$. That is, $f(x)=\sum_{i=0}^{\infty}(f)_ix^i$.
In \cite{Somu1}, the function $H(r,n)=max\{|(d)_r|:d(x)|x^n-1\}$ has been studied. In that paper, we worked on the maximal order of $H(r,n)$ and we showed that \[ H(r,n)\leq (1+o(1))n^{\frac{r(\log 2+o(1))}{\log \log n}}.\]
In this paper, we work on the average order of $H(r,n)$ for a fixed $r$. We will prove the following theorem.
\begin{Main}
	Let $r$ be any natural number then \[ \sum_{n\leq x}H(r,n)\sim c(r)x(\log x)^{2^r-1} \]as $x\rightarrow \infty$ where \[c(r)=.\frac{1}{(2^r-1)!2^rr!}\prod_{p\hspace{1 mm} prime}\bigg(\bigg(1+\frac{2^r-1}{p}\bigg)\bigg(1-\frac{1}{p}\bigg)^{2^r-1}\bigg) .\]
\end{Main}
 \section{Notation}
 Let $\nu(n)$,$\tau(n)$ denote number of prime factors of $n$ and number of divisors of $n$ respectively. Let $\delta(n)$ denote the function $\delta (n)= \left\{\begin{matrix}
 	-1 & n=1 \\ 
 	1 & n>1
 \end{matrix}\right.$. Note that,
\begin{equation}
 \delta(n)\Phi_n(x)=\prod_{d|n}(1-x^d)^{\mu(\frac{n}{d})} .
\end{equation} 
 We say that a formal power series $f(x)=\sum_{n=1}^{\infty}(f)_nx^n$ is dominated by another power series $g(x)=\sum_{n=1}^{\infty}(g)_nx^n$ if $|(f)_n|\leq (g)_n$ for all $n\in \mathbb{N}$.
 If $f(x)$ and $g(x)$ are two formal power series then we denote $f(x)\equiv g(x) \mod x^{r+1}$ if the coefficients of $x^i$ in the power series of $f(x)$ and $g(x)$ are equal for $0\leq i \leq r$.

 \section{Proof of the Main Theorem} 
We require several lemmas in order to prove the main theorem.

\begin{Lemma 1}
	If $f_1(x),\cdots, f_i(x)$ are dominated by $g_1(x),\cdots,g_i(x)$ then $f(x)=\prod_{j=1}^{i}f_j(x)$ is dominated by $g(x)=\prod_{j=1}^{i}g_j(x)$.
\end{Lemma 1}
\begin{proof}
	For all $j\in \{0\}\cup \mathbb{N}$, we have 
	\begin{align*}
	|(f_1\cdots f_i)_j|=|\sum_{h_1+\cdots +h_i=j}(f_1)_{h_1}\cdots (f_i)_{h_i}|.
	\end{align*}
	From triangle inequality and the fact that $f_k(x)$ are dominated by $g_k(x)$ we have 
	\begin{align*}
	|(f_1\cdots f_i)_j|&\leq \sum_{h_1+\cdots+h_i=j}|(f_1)_{h_1}|\cdots|(f_i)_{h_i}|\\
	&\leq \sum_{h_1+\cdots+h_i=j}(g_1)_{h_1}\cdots(g_i)_{h_i}\\
	&=(g_1\cdots g_i)_j.
	\end{align*}
\end{proof}
\begin{Lemma 2}
	There exists a constant $c_1(r)$ depending only on $r$ such that \[ H(r,n)\leq \frac{1}{2^r r!}2^{r\nu(n)}+c_1(r)2^{(r-1)\nu(n)} \] for all natural numbers $n>1$.
\end{Lemma 2}
\begin{proof}
	The proof is similar to that of the proof of Theorem 4.1 of \cite{Somu1}. Let $d(x)$ be a divisor of $x^n-1$. As \[max\{|(d)_r|: d(x)|x^n-1,d(x)\in \mathbb{Z}[x]\}=max\{|(d)_r|: d(x)|x^n-1,d(x)\in \mathbb{Z}[x],d(0)=1\} \] without loss of generality we can assume $d(0)=1$.
	Every divisor $d(x)$ of $x^n-1$ with $d(0)=1$ will be of the form $\prod_{m\in S}\delta(m)\Phi_m(x)$ where $S$ is a subset of set of divisors of $n$.
	Now, \begin{align*}
	d(x)&=\prod_{m\in S}\prod_{d|m}(1-x^d)^{\mu(\frac{m}{d})}\\
	&\equiv \prod_{m\in S}\prod_{\substack{d|m\\d\leq r }}(1-x^d)^{\mu(\frac{m}{d})}\mod x^{r+1}\\
	&\equiv \prod_{d\leq r}(1-x^d)^{\Sigma_1(d)-\Sigma_2(d)}\mod x^{r+1}
	\end{align*}
	where \[\Sigma_1(d)=\sum_{\substack{\mu(\frac{m}{d})=1\\ m\equiv 0 \mod d\\ m\in S}}1\] and \[\Sigma_2(d)=\sum_{\substack{\mu(\frac{m}{d})=-1\\ m\equiv 0 \mod d\\ m\in S}}1.\]
	As $\sum_1(d)=\sum_{\substack{\mu(k)=1\\kd \in S}}1\leq \sum_{\substack{\mu(k)=1\\k|n}}1=2^{\nu(n)-1}$ and $\sum_2(d)=\sum_{\substack{\mu(k)=-1\\kd \in S}}1\leq \sum_{\substack{\mu(k)=-1\\k|n}}1=2^{\nu(n)-1}$, we can conclude that $|\sum_1(d)-\sum_2(d)|\leq 2^{\nu(n)-1}$. Therefore, $(1-x^d)^{\sum_1(d)-\sum_2(d)}$ is dominated by $(1-x^d)^{-(2^{\nu(n)-1})}$. From Lemma 3.1, it follows that 
	\begin{align*}
     |(d(x))_r|&=|(\prod_{d\leq r}(1-x^d)^{\sum_1(d)-\sum_2(d)} )_r|\\
     &\leq \bigg(\prod_{d\leq r}(1-x^d)^{-2^{\nu(n)-1}}\bigg)_r\\
     &=\sum_{c_1+2c_2+\cdots+rc_r=r,c_i\geq 0}	\binom{2^{\nu(n)-1}+c_1}{c_1}\cdots\binom{2^{\nu(n)-1}+c_r}{c_r}.
	\end{align*}
	As \begin{align*}
\binom{2^{\nu(n)-1}+c_1}{c_1}\cdots\binom{2^{\nu(n)-1}+c_r}{c_r}&=\frac{2^{\nu(n)(c_1+\cdots+c_r)}}{2^{c_1+\cdots+c_r}c_1!\cdots c_r!}+O(2^{\nu(n)(c_1+\cdots+c_r-1)})\\ 
&=O(2^{\nu(n)(c_1+\cdots+c_r)})
	\end{align*}
	 and $c_1+2c_2+\cdots+rc_r=r,c_i\geq 0,c_1\neq r$ implies $c_1+\cdots+c_r\leq r-1$ we can conclude that 
	\begin{align*}
	|(d(x))_r|&\leq \binom{2^{\nu(n)-1}+r}{r}+\sum_{c_1+2c_2+\cdots+rc_r=r,c_i\geq 0,c_1\neq r}	\binom{2^{\nu(n)-1}+c_1}{c_1}\cdots\binom{2^{\nu(n)-1}+c_r}{c_r}\\
	&=\frac{2^{\nu(n)r}}{2^rr!}+O(2^{\nu(n)(r-1)}).
	\end{align*}
	Hence there exists a constant $c_1(r)$ depending only on $r$ such that 
	\[ |(d(x))_r|\leq  \frac{1}{2^r r!}2^{r\nu(n)}+c_1(r)2^{(r-1)\nu(n)} .\] Therefore, 	\[ H(r,n)\leq  \frac{1}{2^r r!}2^{r\nu(n)}+c_1(r)2^{(r-1)\nu(n)} .\]
\end{proof}	
\begin{Lemma 3}
There exists a constant $c_2(r)$ depending only on $r$ such that \[H(r,n)\geq \frac{1}{2^r r!}2^{\nu(n)}-c_2(r)2^{(r-1)\nu(n)}\]for all natural numbers $n>1$.	
\end{Lemma 3}
\begin{proof}
	Consider the following divisor $d(x)$ of $x^n-1$, \begin{align*}
	d(x)&=\prod_{\substack{\mu(m)=-1\\m|n}}\delta(m)\Phi_m(x)\\
	&=\prod_{\substack{\mu(m)=-1\\ m|n}}\prod_{d|m}(1-x^d)^{\mu(\frac{m}{d})}\\
	&\equiv \prod_{d\leq r}(1-x^d)^{k(d)}\mod x^{r+1},	
	\end{align*}
	where $k(d)=\sum_{\substack{m|n,\mu(m)=-1\\m\equiv 0 \mod d}}\mu(\frac{m}{d})$. 
	
	Note that, $k(1)=-2^{\nu(n)-1}$ and $|k(d)|\leq 2^{\nu(n)-1}$. Therefore, 
	\begin{equation}
	|(1-x^d)^{k(d)}_{dc_d}|\leq \binom{2^{\nu(n)-1}+c_d}{c_d}.
	\end{equation} 
	We have, 
	\begin{align*}
	(d(x))_r&=(\prod_{d\leq r}(1-x^d)^{k(d)})_r\\
	&=\sum_{c_1+2c_2+\cdots+ rc_r=r,c_i\geq 0}\prod_{d=1}^{r}\bigg((1-x^d)^{k(d)} \bigg)_{dc_d}\\
	&\geq \binom{2^{\nu(n)-1}+r}{r}-\sum_{c_1\neq r,c_1+2c_2+\cdots+rc_r=r}|\prod_{d=1}^{r}\bigg((1-x^d)^{k(d)} \bigg)_{dc_d}|.
	\end{align*}
	From (2), we can conclude that 
	\begin{align*}
	(d(x))_r&\geq  \binom{2^{\nu(n)-1}+r}{r}-\sum_{c_1\neq r,c_1+2c_2+\cdots+rc_r=r}\prod_{d=1}^{r}\binom{2^{\nu(n)-1}+c_d}{c_d}\\
	&=\frac{2^{r\nu(n)}}{2^rr!}+O(2^{(r-1)\nu(n)})
	\end{align*}
	Hence there exists a constant $c_2(r)$ such that 
	$H(r,n)\geq (d(x))_r\geq \frac{2^{r\nu(n)}}{2^rr!}-c_2(r)2^{(r-1)\nu(n)}.$
\end{proof}
\begin{Lemma 4}(See \cite{Murthy}, Problem 4.4.17)
	Let $f(s)=\sum_{n=1}^{\infty}\frac{a_n}{n^s}$, with $a_n=O(n^{\epsilon})$. Suppose that $f(s)=\zeta(s)^kg(s)$, where $k$ is a natural number and $g(s)$ is a Dirichlet series absolutely convergent in $Re(s)>1-\delta$ for some $0<\delta <1$. We have \[\sum_{n\leq x}a_n \sim \frac{g(1)x(\log x)^{k-1}}{(k-1)!} \]
	as $x\rightarrow \infty$.
\end{Lemma 4}
\begin{proof}
	See page 301 of \cite{Murthy} for the proof of the lemma. 
\end{proof}

\begin{Lemma 5}
	For $r\in \mathbb{N}$, we have \[\sum_{n\leq x}2^{\nu(n)r}\sim \prod_{p\hspace{1 mm}prime}\bigg(\bigg( 1+\frac{2^r-1}{p}\bigg)\bigg(1-\frac{1}{p}\bigg)^{2^r-1}\bigg)\frac{x (\log x)^{2^r-1}}{(2^r-1)!} .\]
	\end{Lemma 5}
\begin{proof}
	Let \[f(s)=\sum_{n=1}^{\infty}\frac{2^{r\nu(n)}}{n^s}=\prod_{p \hspace{1 mm} prime}\bigg(\frac{1+\frac{(2^r-1)}{p^s}}{1-\frac{1}{p^s}}\bigg).\]
	As $2^{r\nu(n)}\leq (\tau(n))^r=O(n^{\epsilon})$ we have $2^{r\nu(n)}=O(n^{\epsilon}).$
	Let
	\begin{align*}
	g(s)&=\prod_{p \hspace{1 mm} prime}\bigg(\left(1+\frac{2^r-1}{p^s}\right)\left(1-\frac{1}{p^s}\right)^{2^r-1}\bigg)\\
	&= \prod_{p \hspace{1 mm} prime}\left(1+\frac{d_2}{p^{2s}}+\frac{d_3}{p^{3s}}+\cdots +\frac{d_{2^r-1}}{p^{(2^r-1)s}}\right)
	\end{align*} 
	 for some constants $d_2,d_3,\cdots,d_{2^r-1}$. Notice that Dirichlet's series of $g(s)$ converges absolutely for $Re(s)>\frac{1}{2}$. Observe that $f(s)=\zeta(s)^{2^r}g(s)$. Applying Lemma 3.4, we have \[ \sum_{n\leq x}2^{\nu(n)r}\sim g(1)\frac{x(\log x)^{2^r-1}}{(2^r-1)!}  \] which completes the proof of the lemma.
\end{proof}
Now, we are ready to prove our main theorem.

\begin{proof}
	From lemmas 3.2 and 3.3 we have,
	\begin{equation*}
	 \frac{1}{2^rr!}\sum_{1<n\leq x}2^{\nu(n)r}  -c_2(r)\sum_{1<n\leq x}2^{\nu(n)(r-1)}    \leq \sum_{1<n\leq x}H(r,n) \leq     \frac{1}{2^rr!}\sum_{1<n\leq x}2^{\nu(n)r}   +c_1(r)\sum_{1<n\leq x}2^{\nu(n)(r-1)}. 
	 \end{equation*}
	
	Hence from Lemma 3.5, we get \[\sum_{n\leq x}H(r,n)\sim c(r)x(\log x)^{2^r-1},\] where \[c(r)=.\frac{1}{(2^r-1)!2^rr!}\prod_{p\hspace{1 mm} prime}\bigg(\bigg(1+\frac{2^r-1}{p}\bigg)\bigg(1-\frac{1}{p}\bigg)^{2^r-1}\bigg) .\]
	
\end{proof}


\begin{thebibliography}{99}
	 \bibitem{Bachman}G. Bachman, {\it On the coefficients of cyclotomic polynomials}, Mem. Amer. Math. Soc.,
	 106, no. 510 (1993).
	 
	 \bibitem{Bateman} P. T. Bateman, C. Pomerance, and R. C. Vaughan, {\it On the size of the coeffi-
	 	cients of the cyclotomic polynomial}, Topics in classical number theory, Vol. I,
	 II (Budapest, 1981), Colloq. Math. Soc. J´anos Bolyai, vol. 34, North-Holland,
	 Amsterdam, 1984, pp. 171–202.
	 
	  \bibitem{Bloom}D. M. Bloom, {\it On the coefficients of the cyclotomic polynomials}, Amer. Math. Monthly, 75, 372-377 (1968).
	
	\bibitem{Murthy} M.R. Murthy, {\it Problems in Analytic Number Theory}, Graduate Texts in Mathematics, Springer, 2001.
	
	 \bibitem{Pomerance}C. Pomerance, N.C. Ryan, {\it Maximal height of divisors of $x^n-1$}, Illinois J. Math. {\bf 51} (2007) 597-604.
	
	\bibitem{Somu1}S.T.Somu, {\it On the coefficients of divisors of $x^n-1$}, Journal of Number Theory, Volume {\bf 167}, October 2016, 284-293.
	
	\bibitem{Somu2}  S.T.Somu, {\it On the distribution of numbers related to the coefficients of divisors of $x^n-1$}, Journal of Number Theory, Volume {\bf 170}, January 2017, 3-9.
	
	
\end{thebibliography}
\end{document}